\documentclass[12pt,draft]{amsart}
\usepackage{amssymb, amsmath}
\usepackage{url}

\newtheorem{theorem}{Theorem}[section]
\newtheorem{corollary}[theorem]{Corollary}
\newtheorem{prop}[theorem]{Proposition}
\newtheorem{lemma}[theorem]{Lemma}

\theoremstyle{remark}
\newtheorem{remark}[theorem]{Remark}

\theoremstyle{definition}

\newtheorem{example}[theorem]{Example}

\numberwithin{equation}{section}
\numberwithin{theorem}{section}

\let\abs=\envert

\newcommand{\Lpp}{L{\!}'^{\,p}}
\newcommand{\hpp}{h'^{\,p}}

\newcommand{\intinf}{\int^\infty_{-\infty}}

\newcommand{\R}{{\mathbb R}}

\newcommand{\Sc}{{\mathcal S}}
\newcommand{\fn}{\!:\!}

\providecommand{\abs}[1]{\lvert#1\rvert}
\providecommand{\norm}[1]{\lVert#1\rVert}
\providecommand{\Lany}[1]{L{\!}'^{\,#1}}
\newcommand{\sqrtpt}{\sqrt{\pi t}}
\newcommand{\sqrtp}{\sqrt{\pi}\,}
\newcommand{\fastx}{f\ast\Theta_t(x)}
\newcommand{\fast}{f\ast\Theta_t}

\begin{document}


\date{Preprint September 14, 2023.  To appear in {\it Generalized integrals and applications}
(World Scientific; Eds. T. P. Becerra, H. Kalita, A. Croitoru, B. Hazarika)}
\title[Heat equation]{The heat equation with the $L^p$ primitive integral}
\author{Erik Talvila}
\address{Department of Mathematics \& Statistics\\
University of the Fraser Valley\\
Abbotsford, BC Canada V2S 7M8}
\email{Erik.Talvila@ufv.ca}

\begin{abstract}
For each $1\leq p<\infty$ a Banach space of integrable Schwartz distributions is defined by 
taking the distributional derivative of all functions in $L^p({\mathbb R})$.  
Such distributions can be integrated when multiplied by a function that is the integral of a 
function in $L^q({\mathbb R})$, where $q$ is the conjugate exponent of $p$.
The heat equation on the real line is solved in this space of distributions. 
The initial data is taken to be the distributional derivative of an $L^p({\mathbb R})$ function. 
The solutions are shown to be smooth functions. Initial conditions are taken on in norm. 
Sharp estimates of solutions are obtained and a uniqueness theorem is proved.
\end{abstract}

\maketitle

\section{Introduction}\label{sectionintroduction}

In this paper we study solutions of the heat equation on the real line
with initial value data that is the distributional derivative of an
$L^p$ function
($1\leq p<\infty$).   The solution
is given by convolution of the Gauss--Weierstrass heat kernel with the initial data.  This
produces smooth solutions to the heat equation.   We introduce a norm that is isometric to the $L^p$ 
norm and show
that solutions tend to the initial condition within this norm as $t\to 0^+$.
Since the initial data is the distributional derivative
of an $L^p$ function it can fail to have pointwise values everywhere, can be zero almost everywhere, and
can be the difference of translated Dirac distributions.  And yet, the initial conditions are taken on
in the strong (norm) sense rather than in the weak (distributional) sense. 

For $u\fn\R\times(0,\infty)\to\R$ write $u_t(x)=u(x,t)$.

The classical problem of the heat equation on the real line is, given a function
$F\in L^p$ for some $1\leq p\leq\infty$,  find a function
$u\fn \R\times (0,\infty)\to\R$ such that
$u_t\in C^2(\R)$ for each $t>0$, $u(x,\cdot)\in C^1((0,\infty))$ for each $x\in\R$ and
\begin{align}
&\frac{\partial^2u(x,t)}{\partial x^2}-\frac{\partial u(x,t)}{\partial t}=0 
\text{ for each } (x,t)\in\R\times(0,\infty)\label{heatpde}\\
&\lim_{t\to 0^+}\norm{u_t-F}_p=0.\label{Lpic}
\end{align}
If $p=\infty$ then $F$ is also assumed to be continuous.  Under suitable growth conditions on $u$
the unique solution is given by $u(x,t)=f\ast\Theta_t(x)$, where
the Gauss--Weierstrass heat kernel
is $\Theta_t(x)=\exp(-x^2/(4t))/(2\sqrtpt)$.

Instead of the initial data being an $L^p$ function we take it to be the distributional
derivative of such a function.  A well-defined integration process for these functions
and distributions appears in \cite{talvilaLp}.
The space of derivatives of $L^p$
functions is denoted $\Lpp$.  It is a Banach space isometrically isomorphic to $L^p$.
In $\Lpp$ many properties of $L^p$ functions continue to hold, such as convolution with functions in Lebesgue
spaces, a type of H\"older inequality, continuity of translations, etc.  
Some of the distributions in $\Lpp$ can be identified with signed measures, some can be identified with
smooth functions and some have no pointwise values.

See \cite{cannon} and \cite{widderbook} for classical results on the heat equation, including 
extensive bibliographies.  For initial data in $L^p$ spaces see
\cite{gustafson} and
\cite{hirschmanwidder}.
Distributional solutions
of the heat equation
are considered in \cite{dautraylions} and  \cite{szmydt}, where the
initial data is taken on in the distributional
sense. 
For initial data in the Sobolev space $H^s$,
see \cite{iorio}.  Another paper studying the heat equation
with distributions in Banach spaces is \cite{talvilaheatac}, in which the initial data
is taken to be the distributional derivative of a continuous function.

For $1\leq p<\infty$,
the Lebesgue space on the real line is
$L^p$, which is the set of measurable functions $f\fn\R\to\R$
such that $\intinf |f(x)|^p\,dx<\infty$.  To distinguish between other
types of integrals introduced later, Lebesgue integrals will always
explicitly show the
integration variable and differential as above.

The Schwartz space (test functions), of rapidly decreasing smooth functions, is
denoted $\Sc$ and the tempered distributions $\Sc'$.  See, for example,
\cite{folland} or \cite{friedlanderjoshi}.

The outline of the paper is as follows.

The heat equation is considered in $\Lpp$ in Section~\ref{sectionheatsoln}.  With initial data
$f\in\Lpp$ it is proved that convolution with $f$ and the Gauss--Weierstrass
heat kernel
is a smooth function in $\Lany{s}$ for each $p\leq s<\infty$ but
need not be in $\Lany{s}$ if $1\leq s<p$ (Theorem~\ref{theoremheatsolnproperties}).  
Sharp estimates are given in $\norm{\cdot}'_r$, with $1/p+1/q=1+1/r$.
This convolution gives the unique solution to the heat equation (Corollary~\ref{corollaryheatLppde}).
Existence is proven using a type of Hardy space.
In Proposition~\ref{propcty} continuity is proved with respect to the initial conditions.  

In Section~\ref{sectionestimates} various sharp $L^p$ and pointwise estimates are given.
It is shown that these solutions cannot be of one sign, have limit $0$ at infinity, and have
a vanishing integral over the real line.

An appendix contains a lemma on pointwise differentiation
of convolutions and proof of
some sharp estimates for solutions
of \eqref{heatpde} and \eqref{Lpic}.

\section{Solutions of the heat equation in $\Lpp$}\label{sectionheatsoln}

For $1\leq p<\infty$ let $\Lpp=\{f\in\Sc'\mid f=F'
\text{ for some } F\in L^p\}$.  As explained in \cite{talvilaLp}, if $f\in\Lpp$
then it has a unique primitive in $L^p$.  In this case, $\norm{f}'_p=\norm{F}_p$
and this makes $\Lpp$ and $L^p$ isometrically isomorphic Banach spaces.
As in \cite{talvilaLp}, define
$I^p=\{G\fn\R\to\R\mid
G(x)=\int_0^xg(t)\,dt \text{ for some } g\in L^p\}$

If $f\in\Lpp$ then $f\ast\Theta_t$ is a smooth solution of the heat equation.
The initial value $f$ is taken on in the $\norm{\cdot}'_p$ norm and norm boundedness
ensures uniqueness.  Solutions are shown to be in $\Lany{r}$ for each $p\leq r<\infty$
and sharp norm estimates are given using the results of Theorem~\ref{theoremLpestimates}
in the Appendix.
If a solution $v_t$ of the heat equation has bounded norms $\norm{v_t}'_p$ for $t>0$ then
there exists a distribution $f\in\Lpp$ such that $v_t=\Theta_t\ast f$.  This is proved
using an analogue of a harmonic Hardy space. 
There is continuity with respect to initial conditions.  The example of initial data that
is the difference to two Dirac measures is considered.

\begin{theorem}\label{theoremheatsolnproperties}
Let $1\leq p<\infty$.  Let $f\in\Lpp$ with primitive $F\in L^p$.
Let $t>0$ and define $v(x,t)=v_t(x)=f\ast\Theta_t(x)$.\\
  (a) For each $x\in\R$ we have
$$
v_t(x)  =  \intinf f(y)\Theta_t(x-y)\,dy=F\ast\Theta_t'(x)
=\Theta_t'\ast F(x)=(F\ast\Theta_t)'(x).
$$
The function $v_t$ is uniformly continuous on $\R$.\\
(b) For each $n\geq 1$ we have $v_t^{(n)}(x)=f\ast\Theta_t^{(n)}(x)$ and 
$\partial^n v(x,t)/\partial t^n=f\ast\partial^n\Theta(x,t)/\partial t^n$.
Hence, $v\in C^\infty(\R)\times C^\infty((0,\infty))$.\\
(c) If $p\leq s<\infty$ then $\fast\in\Lany{s}$.\\ 
(d) Let $q,r\in[1,\infty]$
such that $1/p+1/q=1+1/r$.  
There is a constant $K_{p,q}$ such
that $\norm{\fast}'_r\leq K_{p,q}\norm{f}'_p \,t^{-(1-1/q)/2}$ for all $t>0$.  The estimate
is sharp in the sense that
if $\psi\fn(0,\infty)\to(0,\infty)$ such that $\psi(t)=o(t^{-(1-1/q)/2})$
as $t\to 0^+$ or $t\to\infty$ then there is $g\in \Lpp$ such that $\norm{g\ast\Theta_t}'_r/\psi(t)$ is
not bounded.  The constant $K_{p,q}$ cannot be replaced by any smaller number.  In particular,
$\norm{f\ast\Theta_t}'_p\leq\norm{f}'_p$ is a sharp estimate.\\
(e) If $1\leq s<p$ then
$f\ast\Theta_t$ need not be in $\Lany{s}$.
\end{theorem}
Note that if $1/p+1/q=1+1/r$ then, for suitable $q$, 
$r$ can take on any of the values in  $[p,\infty]$,  
with $r=p$ corresponding to $q=1$ and $r=\infty$ corresponding to $p$ and $q$
being conjugates.
And, for suitable $r$, $q$ can take on any of the values in  $[1,p/(p-1)]$, 
with $q=1$ corresponding to
$r=p$ and $q=p/(p-1)$ corresponding to $r=\infty$. If $p=1$ then $q$ can take on all values
in $[1,\infty]$.

The constants $K_{p,q}$ in Theorem~\ref{theoremLpestimates}(b) and part (d) above are the same.
When $r=p$ and $q=1$ the inequality in part (d) reads $\norm{f\ast\Theta_t}'_p\leq\norm{f}'_p$.

\begin{proof}
(a) Since $\Theta_t$ and its derivatives are in each $L^q$ space, the equalities follow by
Theorem~5.1 in \cite{talvilaLp}, Lemma~\ref{lemmaconvolution} in the Appendix,
and usual properties of convolution for Lebesgue integrals.  To prove uniform continuity,
note that if $q$ is the conjugate of $p$ then
$$
\abs{v_t(x)-v_t(y)}  =  \abs{F\ast[\Theta'_t(x)-\Theta'_t(y)]}\leq 
\norm{F}_p\norm{\Theta'_t(x-\cdot)-\Theta'_t(y-\cdot)}_q
\to 0  \text{ as } y\to x.
$$
The last line uses the H\"older inequality and
continuity in the $L^q$ norm.  (When $q=\infty$ this
requires boundedness and uniform continuity of $\Theta'_t$.)

(b) Differentiation in $x$ follows from Lemma~\ref{lemmaconvolution}.
For differentiation in $t$, use part (a) to write
$$
v_t(x)= F\ast\Theta'_t(x)=-\frac{1}{4\sqrtp t^{3/2}}\intinf F(x-y)ye^{-y^2/(4t)}\,dy.
$$
Then
$$
\frac{\partial v(x,t)}{\partial t}  =  \frac{3}{8\sqrtp t^{5/2}}\intinf F(x-y)ye^{-y^2/(4t)}\,dy
   -\frac{1}{16\sqrtp t^{7/2}}\intinf F(x-y)y^3e^{-y^2/(4t)}\,dy.
$$
Differentiation under the integral sign is justified by dominated convergence since, if $t\leq t_0$,
then
$\abs{y^3e^{-y^2/(4t)}}\leq\abs{y}^3e^{-y^2/(4t_0)}$.

(c), (d)  From Theorems 3.6 and 5.1 in \cite{talvilaLp} we have
$$
\norm{\fast}'_r=\norm{F\ast\Theta'_t}'_r=\norm{(F\ast\Theta_t)'}'_r=\norm{F\ast\Theta_t}_r
\leq C_{p,q}\norm{F}_p\norm{\Theta_t}_q.
$$
The rest of the proof is identical to 
Theorem~\ref{theoremLpestimates}(b).

(e) Since $v_t=(F\ast\Theta_t)'\in\Lany{s}$ if and only if $F\ast\Theta_t\in L^s$, this result
is equivalent to Theorem~\ref{theoremLpestimates}(c).
\end{proof}

To prove uniqueness for the initial value problem we use a classical theorem, which follows from
Theorem~6.1 in \cite{widderbook}. 
\begin{theorem}\label{theoremclassicaluniq1}
Let $u\fn \R\times (0,\infty)\to\R$ such that
$u_t\in C^2(\R)$ for each $t>0$, $u(x,\cdot)\in C^1((0,\infty))$ for each $x\in\R$,
$u \in C(\R\times[0,\infty))$,
$\partial^2u(x,t)/\partial x^2-\partial u(x,t)/\partial t=0$ for $(x,t)\in
\R\times(0,\infty)$, $u$ is bounded, $u(x,0)=f(x)$ for
a bounded continuous function $f\fn\R\to\R$. Then the unique solution is given by
$u_t(x)=\fastx$.
\end{theorem}

\begin{corollary}\label{corollaryheatLppde}
Let $1\leq p<\infty$.\\ (a) Let $f\in\Lpp$.
The unique solution of the problem to find
$v\fn \R\times (0,\infty)\to\R$ such that
\begin{align}
&v_t\in\Lpp \text{ for each } t>0\label{cor1}\\
&v_t\in C^2(\R) \text{ for each } t>0, v(x,\cdot)\in C^1((0,\infty)) \text{ for each } x\in\R\label{Lppsmooth}\\
&\frac{\partial^2v(x,t)}{\partial x^2}-\frac{\partial v(x,t)}{\partial t}=0 
\text{ for each } (x,t)\in\R\times(0,\infty)\label{heatLppde}\\
&\norm{v_t}'_p \text{ is bounded }\label{Lpp}\\
&\lim_{t\to 0^+}\norm{v_t-f}'_p=0\label{Lppic}
\end{align}
is given by $v_t=\fast$.\\
(b) If $f\in \Lpp$ and $v_t=\fast$ then $\lim_{t\to 0^+}\norm{v_t}'_p=\norm{f}'_p$.\\
(c) Let $1<p<\infty$ and let $\hpp$ be the set of functions satisfying
\eqref{cor1}, \eqref{Lppsmooth}, \eqref{heatLppde} 
and \eqref{Lpp}.
If $v\in\hpp$ define $\norm{v}'^h_p=\sup_{t>0}\norm{v_t}'_p$.  Let $T\fn\Lpp\to\hpp$ be given by
$T[f]=\Theta_t\ast f$.  Then $T$ is a linear isometric isomorphism.\\
(d) Let $1<p<\infty$. If \eqref{cor1}, \eqref{Lppsmooth}, \eqref{heatLppde} 
and \eqref{Lpp} hold for some $v\in\hpp$ then there is a 
unique distribution $f\in\Lpp$ such that $v_t=\Theta_t\ast f$ and \eqref{Lppic} holds.
\end{corollary}

The proof uses methods from \cite[Theorem~9.2, p.~195]{hirschmanwidder} and \cite[pp.~115-120]{axler}.
Part (c) implies $\hpp$ is a Banach space.
\begin{proof}
(a)
By part (b) of the theorem we can differentiate under the integral sign and this shows $v_t=\fast$ is a smooth
solution of the heat equation.  By part (c) of the theorem, $v_t\in\Lpp$.  
To see that the initial conditions are taken on in the $\norm{\cdot}'_p$
norm,
note that $\norm{v_t-f}'_p=\norm{F\ast\Theta_t-F}_p$.  Since $\intinf\Theta_t(y)\,dy=1$ we have
$$
\norm{F\ast\Theta_t-F}_p=\left(\intinf\left|\intinf\left[F(x-y)-F(x)\right]\Theta_t(y)\,dy\right|^pdx
\right)^{1/p}.
$$
By the Minkowski inequality for integrals, for
example \cite[p.~194]{folland},
$$
\lim_{t\to 0^+}\norm{F\ast\Theta_t-F}_p  \leq  \lim_{t\to 0^+}\intinf\norm{F(\cdot-y)-F(\cdot)}_p\Theta_t(y)\,dy
  =  \lim_{y\to 0}\norm{F(\cdot-y)-F(\cdot)}_p =0.
$$
The last line follows by the fact that $\Theta_t$ is an approximate identity (delta sequence) and by
continuity in the $L^p$ norm.

To prove the solution is unique, suppose there were two solutions $v$ and $w$.  Let $z=v-w$.
Then $z$ satisfies \eqref{cor1}, \eqref{Lppsmooth}, \eqref{heatLppde} and \eqref{Lpp}.  And,
$$
\norm{z_t}'_p\leq\norm{v_t-f}'_p+\norm{w_t-f}'_p\to 0\text{ as } t\to 0^+.
$$

Now define $K\fn\R\to\R$ by $K(x)=1-\abs{x}$ for $\abs{x}\leq 1$ and $K(x)=0$ for $\abs{x}\geq 1$, and let
$K_h(x)=(1/h)K(x/h)$ for $h>0$.  Then $K_h$ is absolutely continuous with compact support.  
Note that $\intinf K(y)\,dy=\intinf K_h(y)\,dy=1$.  We have
\begin{equation*}
K_h'(x)  = \left\{\begin{array}{cl}
0, & \abs{x}>h\\
-{\rm sgn}(x)/h^2, & \abs{x}<h.
\end{array}
\right.
\end{equation*}
Now define $\psi_h(x,t)=K_h\ast z_t(x)=\int_{-h}^hK_h(y)z_t(x-y)\,dy$.  Since $\partial z_t(x)/\partial t$
and $z_t''(x)$ are continuous on compact intervals in $\R\times(0,\infty)$ we
can use dominated convergence to differentiate under the
integral sign and this shows $\psi_h$ is a solution of the heat equation.  
By usual properties of convolutions,
$\psi_h\in C^2(\R)\times C^1((0,\infty))$.  From the H\"older inequality \cite[Theorem~3.6]{talvilaLp},
$\abs{\psi_h(x,t)}\leq\norm{z_t}'_p\norm{K_h'}_q$, where $q$ is the conjugate of $p$.  And,
$\norm{K_h'}_q=2^{1-1/p}/h^{1+1/p}$.  It follows
that $\psi_h$ is bounded on $\R\times(0,\infty)$ and $\lim_{t\to0^+}\psi_h(x,t)=0$.  Define
$\psi_h(x,0)=0$.  The above estimate now shows $\psi_h\in C(\R\times [0,\infty))$.
By Theorem~\ref{theoremclassicaluniq1}, $\psi_h=0$.  Since $z_t$ is continuous, for $(x,t)\in \R\times
(0,\infty)$ we get
$\lim_{h\to0^+}\psi_h(x,t)=0=\lim_{h\to0^+}K_h\ast z_t(x)=z_t(x)\lim_{h\to0^+}\intinf K_h(y)\,dy=z_t(x)$.

(b) This follows from \eqref{Lppic} and the triangle inequality.

(c) Clearly, $T$ is linear.

Show $T$ is an isometry.  From Theorem~\ref{theoremheatsolnproperties}(d), $\norm{T[f]}'^h_p\leq
\sup_{t>0}\norm{f}'_p=\norm{f}'_p$.  And, if $f\in\Lpp$ then
$$
\norm{f}'_p  \leq  \norm{f- \fast}'_p+\norm{\fast}'_p
  \leq  \norm{f- \fast}'_p+\norm{T[f]}'^h_p.
$$
Letting $t\to 0^+$ and using \eqref{Lppic} shows $\norm{f}'_p\leq\norm{T[f]}'^h_p$.
Hence, $T$ is an isometry and necessarily one-to-one.

Now show $T$ is onto $\hpp$.  Let $v\in\hpp$ and find $f\in\Lpp$ such that $T[f]=v$.  For
$h>0$
use notation from the proof of (a) and define $w_h=K_h\ast v_t$. Let $\delta>0$.  Then the
function $(x,t)\mapsto w_h(x,t+\delta)$ satisfies \eqref{cor1},
\eqref{Lppsmooth}, \eqref{heatLppde} 
and \eqref{Lpp} and is continuous and bounded on $\R\times[0,\infty)$.  By 
Theorem~\ref{theoremclassicaluniq1} we have $w_{h}(x,t+\delta)=\intinf\Theta_t(x-y)w_h(y,\delta)\,dy$.
And, $\norm{K_h\ast v_\delta-v_\delta}'_p\to 0$ as $h\to 0^+$ by 
\cite[Theorem~5.3(g)]{talvilaLp}.
(The second sentence in this part of the theorem should read,
Let $F\in L^1$.)
If $q$ is the conjugate of $p$ (\cite[Theorem~3.3(c)]{talvilaLp}) then ${I^q}^\ast=\Lpp$ so,
by weak$^\ast$ convergence, for each $G\in I^q$ we have
$\intinf (K_h\ast v_\delta-v_\delta)G\to 0$.  In particular, take $G=\Theta_t$
and note Definition~3.5 in \cite{talvilaLp}.  Then, also using
the continuity of $v$,
\begin{equation}
\lim_{h\to 0^+}w_h(x,t+\delta)=v_{t+\delta}(x)=\lim_{h\to 0^+}\intinf\Theta_t(x-y)w_h(y,\delta)\,dy=
\intinf\Theta_t(x-y)v_\delta(y)\,dy.\label{tdelta}
\end{equation}

By definition, $v\in\hpp$ is norm-bounded in ${\Lpp=I^q}^\ast$ where $q$ is the conjugate of $p$.
Every norm-bounded sequence in a separable Banach space contains a weak$^\ast$ convergent subsequence.
Since $I^q$ is separable (since $q<\infty$) there is a sequence $\delta_i$ that decreases to $0$ and
a distribution $f\in\Lpp$ such that $\intinf v_{\delta_i}G\to\intinf fG$ for each $G\in I^q$.  In particular,
take $G=\Theta_t$.  Then $\lim_{i\to\infty}\intinf v_{\delta_i}(y)\Theta_t(x-y)\,dy=f\ast\Theta_t(x)$.
Replacing $\delta$ by $\delta_i$ in \eqref{tdelta} and letting $i\to\infty$ now shows
$v_t=T[f]$.

(d)  This follows from (c).
\end{proof}

The following proposition shows continuity with respect to initial conditions.

\begin{prop}\label{propcty}
Let $1\leq p<\infty$.
Let $f,g\in\Lpp$.\\
(a) Let $q,r\in[1,\infty]$
such that $1/p+1/q=1+1/r$.  
Let $K_{p,q}$ be the constant in the proof of Theorem~\ref{theoremheatsolnproperties}.  Then
$$
\norm{\fast-g\ast\Theta_t}'_r\leq K_{p,q}\norm{f-g}'_pt^{-(1-1/q)/2}.
$$
When $r=p$, $q=1$, this reads $\norm{\fast-g\ast\Theta_t}'_p\leq \norm{f-g}'_r.$\\
(b) Let $v,w\fn\R\times(0,\infty)\to\R$ such that
for each $t>0$ we have $v_t$ and $w_t$ in $\Lpp$.  Let $\epsilon>0$.
Suppose $\norm{v_t-f}'_p\to 0$ and $\norm{w_t-g}'_p\to 0$ as $t\to 0^+$.
If $f,g\in\Lpp$ such that $\norm{f-g}'_p<\epsilon$ then for small enough $t>0$ we have
$\norm{v_t-w_t}'_p<2\epsilon$.
\end{prop}
\begin{proof}
(a) This follows from Theorem~\ref{theoremheatsolnproperties}(d) since
$\norm{\fast-g\ast\Theta_t}'_r=\norm{(F-G)\ast\Theta_t}_r$, where $F$ and $G$ are the respective
primitives of $f,g$ in $L^p$.
(b) This follows from the triangle inequality.
\end{proof}

We remark the initial conditions are also taken on in the weak (distributional) sense.
\begin{remark}\label{remarkweaksoln}
If $f\in\Lpp$ with primitive $F\in L^p$ for some $1\leq p<\infty$ then we can define 
$v(x,t)=v_t(x)=f\ast\Theta_t(x)$.  As shown in Corollary~\ref{corollaryheatLppde}, $v$ satisfies the
heat equation and the initial condition is taken on in $\norm{\cdot}'_p$ as $t\to 0^+$.
Such a solution also takes on the initial condition in the weak distributional sense.  To see this
let $\phi\in\Sc$.  Then
$$
\langle v_t-f,\phi\rangle  =\langle (F\ast\Theta_t-F)',\phi\rangle=-\langle F\ast\Theta_t-F,\phi'\rangle
  =-\intinf\left[F\ast\Theta_t(x)-F(x)\right]\phi'(x)\,dx.
$$
We know $\norm{F\ast\Theta_t-F}_p\to 0$ so $\langle F\ast\Theta_t-F,\psi\rangle\to 0$ for each
$\psi\in L^q$ where $q$ is conjugate to $p$.   In particular, this holds when $\psi=\phi'$.
Therefore, the initial condition is taken on in the weak distributional sense.
\end{remark}
Several distributions in $\Lpp$ are given in \cite{talvilaLp}, Section~4., for which \eqref{Lppic} of 
course holds
but \eqref{Lpic} does not hold in any Lebesgue space.  We now consider the example of the difference
of two translated Dirac distributions.
\begin{example}\label{exampledelta}
{\rm
Denote the Dirac distribution supported at $a\in\R$ by $\delta_a$.
Let $F=\chi_{[a,b]}$  then $F$ is in every $L^p$ space.  Therefore, $F'=\delta_a
-\delta_b$ is in every $\Lpp$ space.  Hence, $\Lpp$ is not a subset of any
of the Lebesgue spaces.  Since step functions are dense in $L^p$, linear combinations
of such differences of Dirac distributions are dense in $\Lpp$.

Now let $f=\delta_{-a}-\delta_a$ for fixed $a>0$.  Then $v_t(x)=f\ast\Theta_t(x)=\Theta_t(x+a)-\Theta_t(x-a)$.
Note that $v$ is a solution of the
heat equation in $\R\times(0,\infty)$.  According to Corollary~\ref{corollaryheatLppde} the initial condition 
is taken on as $\norm{v_t-f}'_p\to 0$ as $t\to 0^+$.  

One might wonder in what other sense the initial
conditions are realised.  We've shown above that they are taken on in the weak distributional sense.
For no $1\leq r\leq\infty$ is $\norm{v_t-f}_r$ defined so \eqref{Lpic} is
meaningless.
As a signed Borel measure, we have 
$d\mu(x)=[\Theta_t(x+a)-
\Theta_t(x-a)]dx -d\delta_{-a}(x)+d\delta_a(x)$.  If the Jordan decomposition of $\mu$ is $\mu=
\mu_+-\mu_-$ then $\mu$ converges to $0$ in variation if $V(\mu)=\mu_+(\R)+\mu_-(\R)\to 0$ as $t\to 0^+$.
For any measurable set $E$ we have
\begin{eqnarray*}
V(\mu) &  \geq &  \mu_+(E)+\mu_-(E) \geq \abs{\mu(E)}\\
 & = & \left|\frac{1}{2\sqrt{\pi t}}\int_E\left(
e^{-(x+a)^2/(4t)}-e^{-(x-a)^2/(4t)}\right)dx-\chi_E(-a)+\chi_E(a)\right|.
\end{eqnarray*}
Take $E=(-\infty,-a)$.  Then
\begin{eqnarray*}
V(\mu) & \geq & \frac{1}{2\sqrt{\pi t}}\int_{-\infty}^{-a}\left(
e^{-(x+a)^2/(4t)}-e^{-(x-a)^2/(4t)}\right)dx\\
 & = & \frac{1}{2\sqrt{\pi t}}\int_0^{2a}
e^{-x^2/(4t)}\,dx =\frac{1}{\sqrtp}\int_0^{a/\sqrt{t}} e^{-y^2}\,dy \to \frac{1}{2} \text{ as } t\to 0^+.
\end{eqnarray*}
Hence, $v_t$ does not converge to $f$ in variation.
}

This also shows $v_t$ does not converge to $f$ in the Alexiewicz norm.  See \cite{talvilaregulated}.
\end{example}

\section{Estimates}\label{sectionestimates}

In the previous section we looked at estimates of $\fast$ in $\Lpp$.  Here we look at estimates
in the spaces $L^r$ where $1/p+1/q=1+1/r$.  If $f\in\Lpp$ then $\fast\in L^r$ and
Young's inequality gives estimates in terms of $\norm{f}'_p$.  It is also shown that $\fastx$
cannot be of one sign,
has limit $0$ as $\abs{x}\to\infty$, and has a vanishing integral over the real line.  Pointwise
estimates are given
when $f$ has compact support. 

Let $1\leq q\leq \infty$.  Then
\begin{equation}
\norm{\Theta'_t}_q=\frac{\delta_q}{t^{(2-1/q)/2}} \text{ where } \delta_q=\left\{\begin{array}{cl}
\frac{1}{\sqrt{\pi}}, & q=1\\
\frac{\Gamma^{1/q}((q+1)/2)}{2^{1-1/q}\sqrtp q^{(1+1/q)/2)}}, 
& 1<q<\infty\\
\frac{1}{2^{3/2}\sqrt{\pi e}}, & q=\infty.
\end{array}
\right.\label{Theta'pnormdelta}
\end{equation}
When $q=1$ the integral is elementary.  When $1<q<\infty$ the calculation
uses the definition of the gamma function, $\Gamma(s)=\int_0^\infty e^{-x} x^{s-1}\,dx$. 
When $q=\infty$ this is a calculus minimisation problem. 

\begin{theorem}\label{theoremestimates}
Let $1\leq p<\infty$.  Let $f\in\Lpp$ with primitive $F\in L^p$.
Let $t>0$ and define $v(x,t)=v_t(x)=f\ast\Theta_t(x)$.\\
(a) Let $q,r\in[1,\infty]$
such that $1/p+1/q=1+1/r$.  
There is a constant $L_{p,q}$ such
that $\norm{f\ast\Theta_t}_r\leq L_{p,q}\norm{f}'_p \,t^{-(2-1/q)/2}$ for all $t>0$. 
The estimate
is sharp in the sense that if $\psi\fn(0,\infty)\to(0,\infty)$ such that $\psi(t)=o(t^{-(2-1/q)/2})$
as $t\to 0^+$ or $t\to\infty$ then there is $g\in \Lpp$ such that $\norm{g\ast\Theta_t}_r/\psi(t)$ is
not bounded as $t\to 0^+$ or $t\to\infty$.\\
(b) Fix $t>0$.  Then $v_t$ cannot be of one sign.\\
(c) For each $t>0$, $\intinf v_t(x)\,dx=0$.\\
(d) For each $t>0$, $\lim_{\abs{x}\to\infty}v_t(x)=0$.\\
(e) Suppose $f$ has compact support in $[-R,R]$.  Then
$v_t(x)=O(\abs{x}^{1/p}e^{-(\abs{x}-R)^2/(4t)})$ for fixed $t>0$ as $\abs{x}\to\infty$.  Or,
$$
\abs{v_t(x)} \leq \left\{\begin{array}{cl}
M_p\norm{f}'_1\abs{x} e^{-(\abs{x}-R)^2/(4t)}t^{-3/2}, & p=1,\quad \abs{x}\geq R+\sqrt{2t}\\
M_p
\norm{f}'_p\abs{x}^{1/p}e^{-x^2/(16t)}t^{-(1/2+1/p)},  &
1<p<\infty,\quad \abs{x}\geq 2R,
\end{array}
\right.
$$
where
$$
M_p=\left\{\begin{array}{cl}
\frac{1}{4\sqrtp}, & p=1\\
\frac{3^{1/p}(p-1)^{1-1/p}}{2^{1+2/p}\sqrtp p^{1-1/p}},  &
1<p<\infty.
\end{array}
\right.
$$
\end{theorem}

In part (a) remarks about possible values of $q$ and $r$ are the same as those following 
Theorem~\ref{theoremheatsolnproperties}.

When $r=p$ and $q=1$ the inequality in part (a) reads $\norm{f\ast\Theta_t}_p\leq\norm{F}_p/\sqrt{\pi t}$.
When $r=\infty$ then $p$ and $q$ are conjugates and 
the inequality in part (a) reads $\norm{f\ast\Theta_t}_\infty\leq\delta_q\norm{F}_pt^{-(1+1/p)/2}$.
When $p=1$ then $q=r$ and the sharp constant in (a) is $L_{p,q}=\delta_q$.  This can be proved
similarly to part (b) of Theorem~\ref{theoremLpestimates} using $F=\Theta^\beta_t$ in the limit
$\beta\to\infty$ (effectively, $f=\delta'$).  It is not known in the other cases if the constant
$L_{p,q}$ is sharp or not.

Part (c) appears for $f\in L^1$ as Exercise~8.1.5 in \cite{epstein}.

Positive solutions of the heat equation are known to be convolutions of the Gauss--Weierstrass
heat kernel with an increasing function.  See \cite{widderbook}.

An example for part (d) is
$F(x)=x^{-1/p}\sin(x)$ for $x>1$ and $F(x)=0$ for $x\leq 1$.  Then
$\intinf F(y)\,dy$ exists. And, $F\in L^s$ if and only if $s> p$.

\begin{proof}
(a) By Young's inequality, 
$$
\norm{v_t}_r=\norm{F\ast\Theta'_t}_r\leq C_{p,q}\norm{F}_p\norm{\Theta_t'}_q=\frac{C_{p,q}\norm{f}'_p\delta_q}
{t^{(2-1/q)/2}}.
$$
The constant $C_{p,q}$ is given in the proof of Theorem~\ref{theoremLpestimates}.
We then take $L_{p,q}=C_{p,q}\delta_q$.

As in the proof of Theorem~\ref{theoremLpestimates}, define
$S_t\fn \Lpp\to L^r$ by $S_t[f](x)=f\ast \Theta_t(x)/\psi(t)$.  The estimate
$\norm{S_t[f]}_r\leq C_{p,q}\delta_q\norm{f}'_pt^{-(2-1/q)/2}/\psi(t)$ shows that, for each $t>0$, $S_t$
is a bounded linear operator.  Let $f_t=\Theta'_t$.  Then
$$
\frac{\norm{S_t[f_t]}_r}{\norm{f_t}'_p}
=\frac{\norm{\Theta'_{2t}}_r}{\psi(t)\norm{\Theta_t}_p}
=\frac{\delta_r}{\alpha_p 2^{(2-1/r)/2}\psi(t)t^{(2-1/q)/2}}.
$$
This is not bounded in the limit $t\to 0^+$.  Hence, $S_t$ is not uniformly bounded.
By the Uniform Bounded Principle it is not pointwise bounded.  Therefore, there is a distribution
$g\in \Lpp$ such that $\norm{g\ast\Theta_t}_r\not=O(\psi(t))$ as $t\to 0^+$.  And, the growth estimate
$\norm{v_t}_r=O(t^{-(2-1/q)/2})$ as $t\to 0^+$ is sharp.  Similarly for sharpness
as $t\to\infty$.

(b) Note that $v_t(x)=f\ast\Theta_t(x)=(F\ast\Theta_t)'(x)$ for each $x\in\R$ by 
Theorem~\ref{theoremheatsolnproperties}(a).  If $v_t\geq 0$ then $F\ast\Theta_t$ is increasing.
But $F\ast\Theta_t\in L^p$ by Theorem~\ref{theoremLpestimates}(a) so $F\ast\Theta_t=0$ and 
hence $v_t=0$.

(c) From Example~\ref{exampledelta} we have that linear combinations of differences of Dirac
distributions are dense in $\Lpp$.  Given $\epsilon>0$ there are $-\infty<b_0<b_1<\ldots<b_N<\infty$
and $a_n\in\R$ such that if $g=\sum_{n=1}^Na_n(\tau_{b_{n-1}}\delta-\tau_{b_{n}}\delta)$ then
$\norm{f-g}'_p<\epsilon$.  Then, by Theorem~\ref{theoremestimates}(a),
$$
\left|\intinf f\ast\Theta_t(x)\,dx-\intinf g\ast\Theta_t(x)\,dx\right|\leq \norm{f-g}'_pt^{-(2-1/q)/2},
$$
where $q$ is conjugate to $p$.  It suffices then to show $\intinf g\ast\Theta_t(x)\,dx=0$.  We
have
\begin{eqnarray*}
\intinf g\ast\Theta_t(x)\,dx & = & \intinf\sum_{n=1}^Na_n(\tau_{b_{n-1}}\delta\ast\Theta_t(x)-
\tau_{b_{n}}\delta\ast\Theta_t(x))\,dx\\
 & = & \sum_{n=1}^Na_n\intinf(\Theta_t(x-b_{n-1})-\Theta_t(x-b_{n}))\,dx=0.
\end{eqnarray*}

(d) Let $R>0$. Since $F\in L^1_{loc}$ and $\Theta'_t$ is bounded, dominated convergence shows
$$
\lim_{\abs{x}\to\infty}\int_{-R}^R F(y)\Theta'_t(x-y)\,dy=\int_{-R}^R F(y)\lim_{\abs{x}\to\infty}
\Theta'_t(x-y)\,dy=0.
$$
Given $\epsilon>0$ now take $R>0$ such that $\norm{F\chi_{(R,\infty)}}_p<\epsilon$.  By the H\"older
inequality,
$$
\left|\int_R^\infty F(y)\Theta'_t(x-y)\,dy\right|\leq\norm{F\chi_{(R,\infty)}}_p\norm{\Theta'_t}_q.
$$
Similarly, for integration over $(-\infty,-R)$.

(d)  From the H\"older inequality, $\abs{v_t(x)}\leq\norm{F}_p\norm{\chi_{(x-R,x+R)}\Theta'_t}_q$.
When $p=1$ and $q=\infty$ we compute that if $\abs{x}\geq R+\sqrt{2t}$ then the maximum of 
$\abs{\Theta'_t}$ occurs at
$\abs{x}-R$.  And, $\Theta'_t(\abs{x}-R)=(\abs{x}-R)e^{-(\abs{x}-R)^2/(4t)}/[4\sqrtp t^{3/2}]$.

Now let $1<p<\infty$ and $q$ be its conjugate.  Suppose $x>2R$.  Then
\begin{eqnarray*}
\norm{\chi_{(x-R,x+R)}\Theta'_t}_q^q & = & \frac{1}{(4\sqrtp t^{3/2})^q}\int_{x-R}^{x+R}y^qe^{-qy^2/(4t)}\,dy\\
 & \leq & \frac{(x+R)^{q-1}2t}{(4\sqrtp t^{3/2})^qq}\left[e^{-q(x-R)^2/(4t)}-e^{-q(x+R)^2/(4t)}\right]\\
 & \leq & \frac{3^{q-1}x^{q-1}te^{-qx^2/(16t)}}{2^{q-2}(4\sqrtp t^{3/2})^qq}.
\end{eqnarray*}
\end{proof}

\section{Appendix}\label{sectionappendix}
The Appendix contains a proof of sharp $L^p$ estimates for the heat equation and a lemma
on differentiation of convolutions.

The solution $u(x,t)$ that is in $C^2(\R)$ with respect to $x$ and in $C^1((0,\infty))$ with respect to $t$
of \eqref{heatpde} and \eqref{Lpic} is given by 
the convolution $u_t(x)=F\ast\Theta_t(x)=\intinf F(x-y)\Theta_t(y)\,dy$ 
where the Gauss--Weierstrass heat kernel
is $\Theta_t(x)=\exp(-x^2/(4t))/(2\sqrtpt)$.  For example, see \cite{follandpde}.

The heat kernel has the following properties.  Let $t>0$ and let $s\not=0$ such
that $1/s+1/t>0$.  Then
\begin{align}
\Theta_t\ast\Theta_s&=\Theta_{t+s}\label{thetaconvolution}\\
\norm{\Theta_t}_q&=\frac{\alpha_q}{t^{(1-1/q)/2}} \text{ where } \alpha_q=\left\{\begin{array}{cl}
1, & q=1\\
\frac{1}{(2\sqrt{\pi})^{1-1/q} \,q^{1/(2q)}}, & 1<q<\infty\\
\frac{1}{2\sqrt{\pi}}, & q=\infty.
\end{array}
\right.\label{Thetaqnormalpha}
\end{align}

The last of these follows from the probability integral $\intinf e^{-x^2}\,dx=\sqrt{\pi}$.

We now present some estimates for $F\ast\Theta_t$ when $F\in L^p$, as they are
used in Theorem~\ref{theoremheatsolnproperties} with data in $\Lpp$.  These estimates
are well-known, for example, \cite{gustafson} when $r=\infty$.  
And, \cite[Proposition~3.1]{iwabuchi}.
However, we have not been
able to find a proof in the literature that the estimates are sharp.
\begin{theorem}\label{theoremLpestimates}
Let $1\leq p\leq\infty$ and $F\in L^p$.\\
(a) If $p\leq s\leq\infty$ then $F\ast\Theta_t\in L^s$.\\ 
(b)
Let $q,r\in[1,\infty]$
such that $1/p+1/q=1+1/r$.  
There is a constant $K_{p,q}$ such
that $\norm{F\ast\Theta_t}_r\leq K_{p,q}\norm{F}_p \,t^{-(1-1/q)/2}$ for all $t>0$.  The estimate
is sharp in the sense that if $\psi\fn(0,\infty)\to(0,\infty)$ such that $\psi(t)=o(t^{-(1-1/q)/2})$
as $t\to 0^+$ or $t\to\infty$ then there is $G\in L^p$ such that $\norm{G\ast\Theta_t}_r/\psi(t)$ is
not bounded as $t\to 0^+$ or $t\to\infty$.  The constant 
$K_{p,q}=(c_pc_q/c_r)^{1/2}\alpha_q$, where $c_p=p^{1/p}/(p')^{1/p'}$ with $p, p'$ being
conjugate exponents.  It cannot be replaced with any smaller
number.\\
(c) If $1\leq s<p$ then
$F\ast\Theta_t$ need not be in $L^s$.
\end{theorem}

When $r=p$ and $q=1$ the inequality in part (b) reads $\norm{F\ast\Theta_t}_p\leq\norm{F}_p$.
When $r=\infty$ then $p$ and $q$ are conjugates and 
the inequality in part (b) reads $\norm{F\ast\Theta_t}_\infty\leq\norm{F}_pt^{-1/(2p)}$.

The condition for sharpness in Young's inequality is that both functions be Gaussians.
This fact is
exploited in the proof of part (b).
See \cite[p.~99]{liebloss}, \cite{beckner} and \cite{toscani}.

\begin{proof}
(a), (b)  Young's inequality gives
\begin{equation}
\norm{F\ast\Theta_t}_r\leq C_{p,q}\norm{F}_p\norm{\Theta_t}_q =\frac{C_{p,q}\norm{F}_p\alpha_q}{
t^{(1-1/q)/2}},\label{normyounginequality}
\end{equation}
where $\alpha_q$ is given in \eqref{Thetaqnormalpha}.  
The sharp constant, given in \cite[p.~99]{liebloss},
is $C_{p,q}=(c_pc_q/c_r)^{1/2}$ where $c_p=p^{1/p}/(p')^{1/p'}$ with $p, p'$ being
conjugate exponents. Note that $c_1=c_{\infty}=1$.  Also, $0<C_{p,q}\leq 1$.   We then
take $K_{p,q}=C_{p,q}\alpha_q$.

To show the estimate $\norm{u_t}_r=O(t^{-(1-1/q)/2})$ is sharp as $t\to 0^+$ and $t\to\infty$, let
$\psi$ be as in the statement of the theorem.  Fix $p\leq r\leq\infty$.  Define the family of linear operators
$S_t\fn L^p\to L^r$ by $S_t[F](x)=F\ast \Theta_t(x)/\psi(t)$.  The estimate
$\norm{S_t[F]}_r\leq K_{p,q}\norm{F}_pt^{-(1-1/q)/2}/\psi(t)$ shows that, for each $t>0$, $S_t$
is a bounded linear operator.  Let $F_t=\Theta_t$.  Then, from \eqref{thetaconvolution} and 
\eqref{Thetaqnormalpha},
$$
\frac{\norm{S_t[F_t]}_r}{\norm{F_t}_p}=
\frac{\norm{\Theta_t\ast\Theta_t}_r}{\psi(t)\norm{\Theta_t}_p}
=\frac{\norm{\Theta_{2t}}_r}{\psi(t)\norm{\Theta_t}_p}
=\frac{\alpha_r}{\alpha_p 2^{(1-1/r)/2}\psi(t)t^{(1-1/q)/2}}.
$$
This is not bounded in the limit $t\to 0^+$.  Hence, $S_t$ is not uniformly bounded.
By the Uniform Bounded Principle it is not pointwise bounded.  Therefore, there is a function
$F\in L^p$ such that $\norm{F\ast\Theta_t}_r\not=O(\psi(t))$ as $t\to 0^+$.  And, the growth estimate
$\norm{F\ast\Theta_t}_r=O(t^{-(1-1/q)/2)})$ as $t\to 0^+$ is sharp.  Similarly for sharpness
as $t\to\infty$.

Now show the constant $K_{p,q}$ cannot be reduced.
A calculation shows we have equality in \eqref{normyounginequality}
when $F=\Theta^\beta_t$ and
 $\beta$ is given
by the equation
\begin{equation}
\frac{\beta^{1-1/q}}{(\beta +1)^{1-1/r}}
  =  
\frac{c_pc_q}{c_r}\left(\frac{\alpha_p\alpha_q}{\alpha_r}\right)^2
  =  
\left(1-\frac{1}{p}\right)^{1-1/p}\left(1-\frac{1}{q}\right)^{1-1/q}
\left(1-\frac{1}{r}\right)^{-(1-1/r)}.\label{youngequality}
\end{equation}

First consider the case $p\not=1$ and $q\not=1$.
Notice that $1-1/r=(1-1/q)+(1-1/p)>1-1/q$.  Let
$g(x)=x^A(x+1)^{-B}$ with $B>A>0$.
Then $g$ is
strictly increasing on $(0,A/(B-A))$ and strictly decreasing for $x> A/(B-A)$ so there is a unique
maximum for $g$ at $A/(B-A)$. Put $A=1-1/q$ and $B=1-1/r$.  Then
$$
g\left(\frac{A}{B-A}\right) = \frac{\beta^{1-1/q}}{(\beta +1)^{1-1/r}}
= \left(1-\frac{1}{p}\right)^{1-1/p}\left(1-\frac{1}{q}\right)^{1-1/q}
\left(1-\frac{1}{r}\right)^{-(1-1/r)}. 
$$
Hence, \eqref{youngequality} has a unique positive solution for $\beta$ given by
$\beta=(1-1/q)/(1-1/p)$.  

If $p=1$ then $q=r$.  In this case, \eqref{youngequality} reduces to
$(1+1/\beta)^{1-1/q}=1$ and the solution is given in the limit $\beta\to\infty$.  Sharpness
of \eqref{normyounginequality} is then given in this limit.  It can also be seen that taking $F$ to be the Dirac
distribution gives equality.

If $q=1$ then $p=r$.  Now, \eqref{youngequality} reduces to $(\beta+1)^{1-1/p}=1$ and $\beta=0$.  There
is equality in \eqref{normyounginequality} when $F=1$.  This must be done in the limit $\beta\to0^+$.

If $p=q=r=1$ then there is equality in \eqref{normyounginequality} for each $\beta>0$.

Hence, the constant in \eqref{normyounginequality} is sharp.

(c) Suppose $F\geq 0$ and $F$ is decreasing on $[c,\infty)$ for some $c\in\R$.  Let $x>c$.  Then
\begin{eqnarray*}
F\ast\Theta_t(x) & \geq & \int_c^x F(y)\Theta_t(x-y)\,dy
 \geq  F(x)\int_c^x \Theta_t(x-y)\,dy\\
 & = & \frac{F(x)}{\sqrt{\pi}}\int_0^{(x-c)/(2\sqrt{t})} e^{-y^2}\,dy
 \sim  F(x)/2 \quad \text{ as } x \to \infty.
\end{eqnarray*}
Now put $F(x)=1/[x^{1/p}\log^2(x)]$ for $x\geq e$ and $F(x)=0$, otherwise.
For $p=\infty$ replace $x^{1/p}$ by $1$.
\end{proof}

We prove a lemma on pointwise differentiation of $L^p$ convolutions.  
A similar result for $p=1$ and functions
with bounded derivatives appears
in \cite[Proposition~8.10]{folland}.

\begin{lemma}[Differentiation of Convolutions]\label{lemmaconvolution}
Let $F\in L^p$ for some $1\leq p\leq \infty$.  Let $q$ be the conjugate of $p$ and suppose
$\psi\in L^q$.\\
(a) Suppose $\psi$ is a differentiable function such that $\psi'$ is absolutely continuous
and $\psi''\in L^q$.
Then $F\ast\psi$ is differentiable and $(F\ast\psi)'(x)=F\ast\psi'(x)$ for each $x\in\R$.\\
(b) Suppose $\psi$ is a differentiable function such that $\psi$ is absolutely continuous, $\psi'$
is bounded on each compact interval and $\psi'\in L^q$.
Then $F\ast\psi$ is differentiable and $(F\ast\psi)'(x)=F\ast\psi'(x)$ for each $x\in\R$.\\
(c) Suppose $\psi$ is a $C^\infty$ function such that $\psi^{(n)}\in L^q$ for each $n\geq 0$.
Then $F\ast\psi\in C^\infty(\R)$ and $(F\ast\psi)^{(n)}(x)=F\ast\psi^{(n)}(x)$ for each $n\geq 1$ 
and each $x\in\R$.
\end{lemma}
\begin{proof}
Let $w=F\ast\psi$. (a)  Without loss of generality, we can take $h>0$ and consider
$$
\frac{w(x+h)-w(x)}{h}  =  \intinf F(y)\left[\frac{\psi(x+h-y)-\psi(x-y)}{h}\right]dy
  =  F\ast\psi'(x)-R(h).
$$
The remainder from Taylor's theorem is
$$
R(h)  =  \frac{1}{h}\intinf F(y)\int_{x-y}^{x+h-y}(z-x-h+y)\psi''(z)\,dz\,dy.
$$
And, by the H\"older inequality, Jensen's inequality and the Fubini--Tonelli theorem,
\begin{eqnarray*}
\abs{R(h)} & \leq & \norm{F}_p\left(\intinf\int_{x-y}^{x+h-y}\abs{z-x-h+y}^q\abs{\psi''(z)}^q\,\frac{dz}{h}
\right)^{1/q}dy\\
 & = & \norm{F}_p\left(\intinf\abs{\psi''(z)}^q\int_{x-z}^{x+h-z}\abs{z-x-h+y}^q
\,\frac{dy}{h}\,dz\right)^{1/q}\\
 & = & \frac{\norm{F}_p\norm{\psi''}_q\abs{h}}{(q+1)^{1/q}}\to 0 \text{ as } h\to 0.
\end{eqnarray*}
Minor changes are needed when $p=1$.
Therefore, $(F\ast\psi)'(x)=F\ast\psi'(x)$.

(b) Let $R>0$.  Write
$$
F\ast\psi(x)=\int_{-R}^R F(y)\psi(x-y)\,dy+\int_{\abs{y}>R}F(y)\psi(x-y)\,dy.
$$
Since $F\in L^1_{loc}$, and $\psi'$ is bounded on compact intervals, dominated convergence shows
$$
\frac{d}{dx}\int_{-R}^R F(y)\psi(x-y)\,dy=\int_{-R}^R F(y)\psi'(x-y)\,dy.
$$
Now, let $h>0$ and consider
\begin{align*}
&\left|\int_R^\infty F(y)\left[\frac{\psi(x+h-y)-\psi(x-y)}{h}\right]dy\right|  \leq  \\
&\qquad
\norm{F\chi_{(R,\infty)}}_p\left(\int_R^\infty \left|\frac{\psi(x+h-y)-\psi(x-y)}{h}
\right|^q dy\right)^{1/q}.
\end{align*}
Given $\epsilon>0$ take $R$ large enough so that $\norm{F\chi_{(R,\infty)}}_p<\epsilon$.
Using the fundamental theorem of calculus, Jensen's inequality and then the Fubini--Tonelli theorem, 
we have
\begin{eqnarray*}
\int_R^\infty \left|\frac{\psi(x+h-y)-\psi(x-y)}{h}
\right|^q\,dy & = & \int_R^\infty\left|\int_{x-y}^{x+h-y}\psi'(z)\,\frac{dz}{h}\right|^q\,dy\\
 & \leq & \int_R^\infty\int_{x-y}^{x+h-y}\left|\psi'(z)\right|^q\,\frac{dz}{h}\,dy\\
 & = & \intinf\left|\psi'(z)\right|^q\int_{x-z}^{x+h-z}\frac{dy}{h}\,dz\\
 & = & \norm{\psi'}_q^q.
\end{eqnarray*}
Therefore,
$$
\left|\int_R^\infty F(y)\frac{\psi(x+h-y)-\psi(x-y)}{h}\,dy\right|  \leq \epsilon\norm{\psi'}_q.
$$
The case $p=1$ is similar.

(c) Follows from (a) or (b).
\end{proof}

\end{document}